\documentclass{amsart}

\usepackage{amsmath}
\usepackage{amssymb}

\def \bbN {{\mathbb N}} 
 
\def \bbZ {{\mathbb Z}}

\def \Sym {{\rm Sym}}
\def \Stab {{\rm Stab}}

\def \Image {{\rm Im}} 
\def \calB {\mathcal{B}} 
\def \Aut {{\rm Aut}} 
\def \ker {{\rm Ker}} 

\newcounter{iii}

\begin{document}

\title{Computing Symmetric Normalisers}

\author{Andreas-Stephan Elsenhans}

\begin{abstract}
The computation of the normaliser of a permutation group
in the full symmetric group is an important and hard problem in computational group theory. 
This article reports on an algorithm that builds a descending chain of overgroups 
to determine the normaliser. A detailed performance test in {\tt magma} 
shows the improvement for large examples of  intransitive and imprimitive groups.
\end{abstract}

\maketitle

\section{Introduction}
Computing intersections, testing subgroups for conjugacy, determining the stabiliser of a set, and
determining the normaliser of a given subgroup are known to be hard problems in computational group theory.
Efficient algorithms are still a research problem. 
An overview of different approaches is given in~\cite{HEB}. 
Most notable is the work of Leon~\cite{Leon}. 

A method to compute the normaliser $N_G(U)$  starts 
with the exact sequence of centraliser, normaliser and automorphism group
$$
0 \rightarrow C_G(U) \rightarrow N_G(U) \rightarrow \Aut(U)\, .
$$
Here, an element $g$ of the normaliser is mapped to the conjugation map $u \mapsto g u g^{-1}$.
The centraliser can be computed with the methods \cite[Sec.~4.6.5]{HEB}.
After a determination of the automorphism group of $U$, one needs to check which of the 
automorphisms can be realised by  conjugation. 
In particular cases such as regular permutation group or almost simple groups, 
this approach is very efficient.

A second approach starts with
$$
N_G(U) =  \!\!\!\!\!\!\!\! \bigcup_{\tiny \begin{array}{c}  gU \in G / U \\ U = gUg^{-1} \end{array}} \!\!\!\!\!\!\!\! g U \, .  
$$
Thus, one could in principle compute $N_G(U)$ by an inspection of all the cosets in $G / U$.
A highly optimised backtrack search turns this into a practical algorithm. 
Therefore, one could use the index $[G : U]$ as a crude upper bound for the run time.

As current algorithms are randomised, repeated computations with the same input data may result in 
very different timings. Thus, it is already hard to determine the performance in 
practice. Not to mention theoretical considerations.
But, one can roughly say that normaliser algorithms perform well for 
primitive permutation groups. The slowest examples are intransitive groups. 
The performance for transitive but imprimitive groups is in between. 

The above indicates that one should start 
with a close inspection of intransitive permutation groups and build a specialised 
normaliser algorithm for them. Furthermore, strategies successfully applied to intransitive 
groups should be carried over to imprimitive groups.
Very recently, work in this direction was carried out by Mun See Chang~\cite{Ch}. 
The reader interested in normalisers of primitive groups of very large degree 
may consult~\cite{RS}.

A second source of input is the authors experience in Galois group computations 
and the methods to analyse permutation groups for that.
The general strategy for this is to build a chain of smaller and smaller 
overgroups, until the actual Galois group is reached. 
An indication that this strategy might help for normalisers is the following. 
Assume that one wants to find $N_{S_n}(G)$ for $n \geq 5$. Then only the special cases 
$G = S_n, A_n, \{ e \}$, result in the full symmetric group $S_n$. In all other cases,
the result will be a proper subgroup and not $A_n$. These subgroups are either intransitive, 
imprimitive, or primitive and of huge index in $S_n$~\cite[Thm.~5.6.C]{DM}. 
Thus, if one computed a non-trivial overgroup of the normaliser, one could try to use the 
structure of its orbits or its blocks to further shrink the overgroup. 
In the case of a primitive overgroup, one benefits 
from the fact that its order is magnitudes smaller than the order of~$S_n$. 
The smaller the final overgroup is, the faster the normaliser can be
found by a backtrack search.

The article is structured as follows. First, intransitive groups are inspected. 
This leads to a non-trivial overgroup of the normaliser by only carrying out 
computations with the orbit images. This is complemented by a normaliser 
overgroup algorithm for imprimitive transitive groups. It determines an 
overgroup via the automorphism group of a graph that encodes
the structure of the block systems. This automorphism group is found quickly 
by calling {\tt nauty}~\cite{MP}.

Next, methods to shrink the overgroup are presented. The general strategy is to map the 
groups to more accessible quotients. For this, a list of maps based on block systems is given. 
In the case of an intransitive group $G$, all maps $\phi$ defined on the given 
overgroup are classified of the kind that $\phi(G)$ has only orbits of prime length and the orbit 
groups are 1-dimensional affine groups. In this particular case, the normaliser relates to 
the automorphism group of a linear code. For the latter, efficient algorithms are available 
that can accelerate the normaliser computation. 

Finally, a systematic performance test of a {\tt magma} implementation was carried out.
It shows that knowing the constructed overgroup of the normaliser speeds up
the determination of the exact normaliser. The code will be available in the
next version of {\tt magma}.

\section{An approach for intransitive groups}

\subsection{Notation and Terminology}
In this article, the following notation and terminology will be used:
\begin{enumerate}
\item
Given a permutation group $G \subset \Sym(M)$, the set $M$ is called the {\em domain} of $G$.
The {\em degree} of a permutation group is the size of its domain.
\item
The normaliser of $G$ in $\Sym(M)$ is called the 
{\em symmetric normaliser} of $G$.
\item
Let $\Delta$ be an invariant subset of the domain of a permutation group $G$. Then  
$\pi_\Delta \colon G \rightarrow \Sym(\Delta)$ denotes
the projection of $G$ onto its action on $\Delta$.
\item
Two permutation groups are called {\em permutation isomorphic} if a bijection of their domains exists
that induces an isomorphism of the groups.
\item
A subgroup  $U \subset G_1 \times G_2 \times \cdots \times G_n$ is called a {\it subdirect product}
if all projections $\pi_i \colon U \rightarrow G_i$ are surjective.
\item
A direct product of permutation groups is viewed as a intransitive permutation group acting on the
disjoint union of the domains of the factors.
\item
A wreath product $H \wr G$ of two permutation groups of degrees $d_H$ and $d_G$ is an imprimitive
permutation group of degree $d_H \cdot d_G$~\cite[Sec.~2.6]{DM}.
\item
A {\em linear code} in $(\bbZ / p \bbZ)^n$ is a linear subspace. The automorphism group of a linear code
is the group of all the maps given by monomial matrices that map the code onto itself. 
In the special case of $p=2$, this group consists only of permutation matrices and is 
therefore viewed as a subgroup of $S_n$.
\end{enumerate}

\subsection{Fact}
Let $G$ be a group acting transitively on a set $\Omega$ and let $H \triangleleft G$ be a normal subgroup. Then 
\begin{itemize}
\item[(i)] the orbits of $H$ form a system of blocks for $G$.
\item[(ii)] If $\Delta$ and $\Delta^\prime$ are two $H$-orbits then 
$\pi_\Delta (H)$ and $\pi_{\Delta^\prime} (H)$ are permutation
isomorphic.
\end{itemize}

\proof
This is part of~\cite[Thm 1.6A]{DM}.
\qed

\subsection{Fact\label{WPemb}}
Let $G$ be a transitive permutation group with block system\break $B_1,\ldots,B_k$. 
Then $G$ can be embedded into
$$
\pi_{B_1}(\Stab_G(B_1)) \wr S_k
$$
by a suitable renumbering of its domain.
\proof
This is~\cite[Thm 1.8]{Ca} or~\cite[Satz 1.10]{G}. \qed

\subsection{Application}
The results above give an overgroup of the normaliser of a permutation group.
For this, let $G_1,\ldots,G_m$ be pairwise non permutation isomorphic transitive permutation 
groups of degree $d_1,\ldots,d_m$ and
$$
U \subset G_1^{e_1} \times \cdots \times G_m^{e_m} 
$$
be a subdirect product. If we view $U$ as a permutation group acting on the disjoint union of the
domains of the factor, its symmetric normaliser is contained in
$$
N_{S_{d_1}}(G_1) \wr S_{e_1} \times \cdots \times  N_{S_{d_m}}(G_m) \wr S_{e_m}\, .
$$
As normalisers in smaller degree can be computed recursively, this gives a starting point 
for the normaliser computation of $U$. 

\subsubsection*{Remark} A similar result is obtained in~\cite[Chap.~4]{Ch}.

\subsection{Algorithm\label{OrbitSorting}}(Orbit sorting)\\
{\bf Input:} A permutation group $G$ of degree $n$. \\
{\bf Output:} A relabelling of the domain of $G$ such that $G$ is a subdirect product of distinct transitive factors. 
I.e., $G \subset G_1^{e_1} \times \cdots \times G_m^{e_m}$. \\
{\bf Steps:}
\begin{enumerate}
\item
Compute the orbits $\Delta_1,\ldots,\Delta_k$ of $G$ and the orbit images $\pi_{\Delta_i}(G)$.
\item
For each pair $(i,j) \in \{1,\ldots,k\}^2$ test, if $\pi_{\Delta_i}(G)$ and $\pi_{\Delta_j}(G)$ are permutation isomorphic.
\item
Relabel the domain of $G$ s.t.
\begin{enumerate}
\item
The elements in each orbit are labelled with consecutive integers.
\item
Orbits with permutation isomorphic orbit images immediately follow each other.
\item
If the orbit images on the orbits $\{k,\ldots,k+l-1\}$ and  $\{k+l,\ldots,k+2l-1\}$
are permutation isomorphic, then an isomorphism is given by $k+j \mapsto k+l+j$ for $j =0,\ldots,l-1$.
\end{enumerate}
\item
Return $G$ with the new labelling. Also return the direct product of the orbit images 
as the group that contains $G$ as a subdirect product. 
\end{enumerate}

\subsection{Remark}
In order to implement the algorithm above in a computer algebra system, one has to test permutation isomorphy.
As permutation isomorphic groups have the same degree, this can be done as follows.

Relabel the domains of both groups by $1,\ldots,n$. Test whether this yields two conjugate subgroups of $S_n$.
As the test for conjugacy also gives an element that maps one group to the other, 
a permutation isomorphism of the initial groups can be constructed. If the groups are not conjugate then
they are not permutation isomorphic.

\subsection{Algorithm\label{DPWP}}(Overgroup of the normaliser)\\
{\bf Input:} A permutation group $G$. \\
{\bf Output:} An overgroup of the symmetric normaliser of $G$. \\
{\bf Steps:}
\begin{enumerate}
\item
Use Algorithm~\ref{OrbitSorting} to get a labelling of the domain of $G$ that 
presents $G$ as a subdirect product of 
$$
G_1^{e_1} \times \cdots \times G_m^{e_m}\, ,
$$
for distinct transitive permutation groups $G_1,\ldots,G_m$ of degree $d_1,\ldots,d_m$.
\item
Return 
$$
\prod_{j = 1}^m N_{S_{d_j}}(G_j) \wr S_{e_j}
$$
as an overgroup of the normaliser with respect to the new labelling of the domain.
\end{enumerate}

\subsection{Remark}
Given a direct product of transitive permutation groups, this overgroup 
coincides with the normaliser.

\section{An approach for transitive groups}

\subsection{Introduction}
The approach for intransitive groups described above was dri\-ven by the orbit structure. 
For transitive groups, the block systems have to be used instead. This is limited to
imprimitive groups. But as mentioned in the introduction, 
primitive groups are well covered by other methods.

\subsection{Definition}
Let $G \subset S_n$ be a transitive permutation group with  block system $\calB$.
The block system $\calB$ is said to be {\it generated by the pair} $(i,j) \in \{1,\ldots,n\}^2$, 
if it is the finest block system of $G$ such that $i$ and $j$ are contained in the
same block. A block system that is generated by a single pair is called {\it principal}.

\subsection{Remarks}
\begin{enumerate}
\item
As an algebraic structure, the block systems of $G$ form a lattice. This lattice is generated by the
principal block systems. The two operations of the lattice are given by intersecting and merging the 
block systems.
\item
A transitive group of degree $n$ has at most $n-1$ principal block systems. They are generated by the 
pairs $(1,2),\ldots,(1,n)$. The generating pair may not be unique.
\item
The total number of block systems of a permutation group is not bounded by a polynomial in the degree. 
This can be shown considering the regular representations of~$(\bbZ / 2 \bbZ)^k$, for $k \in \bbN$.
\end{enumerate}

\subsection{Lemma} \label{NormalisatorOperation}
Let $G \triangleleft N \subset S_n$ be transitive permutation groups.
Further, let a block system $(B_1,\ldots,B_k)$ of $G$ and an element $\sigma \in N$ be given. 
For $(C_1,\ldots,C_k) := (\sigma B_1,\ldots,\sigma B_k)$, we have:
\begin{enumerate}
\item 
The system $(C_1,\ldots,C_k)$ is a block system for $G$.
\item
The block system $(B_1,\ldots,B_k)$ is principal if and only if $(C_1,\ldots,C_k)$ is principal.
\item
The groups $\Stab_G(B_1)$ and  $\Stab_G(C_1)$ are permutation isomorphic.
Furthermore, $\pi_{B_1} (\Stab_G(B_1))$ and $\pi_{C_1} (\Stab_G(C_1)))$ are permutation isomorphic, too. 
\item
Denote by $\phi_B$ and $\phi_C$ the homomorphisms $G \rightarrow S_k$, given by the actions on the block systems
$(B_1,\ldots,B_k)$ and $(C_1,\ldots,C_k)$ respectively. 
The images $\phi_B(G)$ and $\phi_C(G)$ are permutation isomorphic, as well.
\end{enumerate}

\proof
As the conjugation by $\sigma$ is a permutation automorphism of $G$, the claim follows. \qed

\subsection{Remarks}
\begin{enumerate}
\item
In the special case that $G$ has exactly one block system with $k$ blocks of size 
$\ell$, the normaliser has this block system as well and is therefore contained in 
the wreath product $S_\ell \wr S_k$.
\item
If the normaliser has several block systems then it is contained in the intersection of 
the corresponding wreath products.
\item
If $G$ has several block systems of equal block size then the normaliser may permute the 
block systems. In particular, the normaliser of the regular permutation representation 
of $(\bbZ / 2 \bbZ)^d$ is the primitive affine group in dimension $d$. 
\end{enumerate}
The next goal is to use the results above to determine an overgroup of the normaliser.
This is done by encoding the block systems in a graph.

\subsection{Definition}
Let $\calB_1 := (B_{1,1},\ldots,B_{1,n_1}),\ldots, \calB_m := (B_{m,1},\ldots,B_{m,n_m})$ be block 
systems partitioning $\{1,\ldots,n\}$. Then the graph described below with $n+n_1 + \cdots + n_m + m$ vertices is 
called the {\em block system graph} of the given block systems.

The set of vertices is the union of the set of point $\{1,\ldots,n\}$, the set of all blocks 
$\{B_{i,j}\}$, and the set of the block systems $\{ \calB_1,\ldots, \calB_m\}$. 
The edges are given by the following rules
\begin{enumerate}
\item
A point vertex $k$ is joined with the block vertex $B_{i,j}$ if and only if $k \in B_{i,j}$.
\item
A block vertex $B$ is joined with a block system vertex $\calB$ if and only if $B$ is a block of $\calB$.
\end{enumerate}
For a transitive permutation group $G$, we call the block system graph constructed from 
all the block systems of $G$ the {\em block system graph of $G$}. Furthermore, the graph constructed 
from all  principal block systems is called the {\it principal block system graph of $G$}.

\subsection{Remark}
For a transitive subgroup $G \subset S_n$, Lemma~\ref{NormalisatorOperation} implies that $N_{S_n}(G)$ acts on the 
block system graph of $G$ and the principal block system graph of~$G$.

Thus, given a transitive, but imprimitive permutation group $G$, the automorphism group of the 
block system graph and the principal block system 
graph are overgroups of the symmetric normaliser of $G$. 
This results in the following algorithm.

\subsection{Algorithm\label{GraphGroup}}(Normaliser overgroup by block system graph)\\
{\bf Input:} A transitive, imprimitive permutation group $G \subset S_n$. \\
{\bf Output:} An overgroup of the symmetric normaliser of $G$. \\
{\bf Steps:}
\begin{enumerate}
\item
Compute all the principal block systems of $G$.
\item
Build the principal block system graph of $G$.
\item
Colour the block system vertices with one colour, 
the block vertices with a second colour, and the point vertices with a third colour.
\item
Compute the group of colour preserving automorphisms of the graph.
\item
Compute its action on the point vertices $\{1,\ldots,n\}$.
\item
Return the image as an overgroup of the normaliser.
\end{enumerate}

\subsection{Remarks}
\begin{enumerate}
\item
One can refine the algorithm above by colouring the block system vertices 
$\calB_1,\ldots,\calB_k$ of the block system graph with 
invariants associated with the block systems. 
For the block system $\calB$ containing the block $B$, one can use 
any property of $\Stab_G(B)$, $\pi_{B}(\Stab_G(B))$, or $\phi_\calB(G)$. 
The authors experiments indicate that
$\#\phi_\calB(G)$ and the structure of the abelian quotient 
of $\Stab_{B}(G)$ are helpful choices.
\item
If we assign a special colour to one point vertex of the block system graph then the result 
will be the  point stabiliser $S_1$ of the overgroup.
But, as $G$ is transitive and contained in its normaliser, the composition of $G$ and 
the computed stabiliser $S_1$ will be the overgroup of the normaliser searched for. 
In some cases, this modified colouring results in a significant speed up of the 
computation of the graph automorphism group. 
\end{enumerate}

\section{Refining the overgroup}

\subsection{Introduction}
Using the techniques above, one can compute an overgroup $S$ of the normaliser, for every 
non-primitive permutation group $G$. At this point, one could call any normaliser algorithm 
and compute the normaliser of $G$ in $S$. 
But, at least when the index of $G$ in $S$ is large, 
it is more efficient to search for smaller intermediate groups first. 

\subsection{Lemma} \label{Lemma_trivial}
Let $V,U, N \subset G$ be subgroups and assume that $N$ is normal in $G$. Then
$U \cap N$ is normal in $U$. Furthermore, if $\phi \colon G \rightarrow H$
is a homomorphism then $\phi(N)$ is normal in $\phi(G)$.
Finally, we have $N_G(V) \subset N_G(V \cap N)$.

\proof 
As the kernel of the projection homomorphism $G \rightarrow G / N$ restricted to $U$ is $U \cap N$, we can conclude that 
$U \cap N$ is normal in $U$. 
Furthermore, for $\phi(g) \in \phi(G)$ we have 
$$
\phi(g) \phi(N) \left(\phi(g)\right)^{-1} = \phi(g N g^{-1}) = \phi(N)\,.
$$
Thus, $\phi(N)$ is normal in $\phi(G)$.
Finally, as $V$ is normal in $N_G(V)$, the first part implies that $V \cap N$ is normal in $N_G(V)$. Thus, 
$N_G(V) \subset N_G(V \cap N)$ follows.
\qed

\subsection{Proposition\label{F-Case}} 
Let $G \subset S_k \wr S_\ell$ be a transitive permutation group. 
Furthermore, let $B = \{1,\ldots,k\}$ be a block.
Then 
$$
N_{S_k \wr S_\ell}(G) \subset N_{S_k}\left(\pi_B (\Stab_G(B)) \right) \wr S_\ell 
$$
is an overgroup of the normaliser of $G$ in $S_k \wr S_\ell$.

\proof
This is a direct consequence of Theorem~\ref{WPemb} and Lemma~\ref{Lemma_trivial}, 
applied to the stabiliser of $B$ and its action on $B$.  
\qed

\subsection{Application}
Once it is established that the normaliser of a transitive group has a block system, the normaliser is contained in 
the wreath product corresponding to the block system. Therefore, Proposition~\ref{F-Case} can be used to refine the
normaliser overgroup by intersecting it with the wreath product.

\subsection{Remark}(Accessible homomorphisms) \\
In general, quotients of permutation groups can be hard to handle. In particular, the quotient may not have a 
permutation representation in small degree. A quotient and the corresponding homomorphism are considered as 
accessible when it can be constructed and used efficiently. Examples for this include
\begin{enumerate}
\item
the action of an imprimitive group on a block system,
\item
the action of an intransitive group on an orbit,
\item
the direct product of the actions on block systems of orbit images.
In particular, one can take one minimal block system of each orbit image.
\item
The action of an intransitive group on its first $k$ orbits.
\end{enumerate}
Furthermore, the following proved to be efficient. 

\subsection{Remark\label{BlockQuotient}}
Let $U \subset G$ be permutation groups. Denote by $\phi$ the permutation representation of $G$
given by the coset action on $G/U$. 
This coset action can be used
to build a representation of the wreath product $G \wr S_n$. 
For this, consider $T := U \times \left( G \wr S_{n-1} \right)$ as a subgroup of $G \wr S_n$.
Then the coset action of $G \wr S_n$ on $(G \wr S_n) / T$ is the induced permutation representation of $\phi$.

In practice, this can be computed as follows. Let $B$ be a block of the permutation group $W$ and 
$U \subset \pi_B(\Stab_W(B))$ a subgroup. Then the action of $W$ on $W / \pi_B^{-1}(U)$ is the induced action
of the action of $\pi_B(\Stab_W(B))$ on $\pi_B(\Stab_W(B)) / U$. 
This is sometimes called the {\it block quotient}.

\subsection{Application}
A surjective homomorphism  $\varphi \colon G \rightarrow H$ can be used to generate a normaliser overgroup via
$$
N_G(U) \subset \varphi^{-1} \left( N_H(\varphi(U)) \right) \, .
$$
If the right hand side can be computed quickly and is a proper subgroup of $G$ then the group obtained 
is a smaller normaliser overgroup.

\section{Using automorphisms of codes\label{CodeAlg}}

\subsection{Introduction}
Given a permutation group with $k$ orbits of size 2, the approach above locates the normaliser in the 
group $S_2 \wr S_k$. Similarly, any subdirect product of $A_3^k$  
results in the normaliser overgroup $S_3 \wr S_n$. 

Thus, the index of the normaliser in the overgroup may still be large. Furthermore, the  overgroup
does not have any helpful quotients and a direct call of a backtrack normaliser algorithm may 
use a lot of time.

But, viewing this type of input group as a general permutation group does not suit its structure. 
The structure of a subdirect product in $S_2^k$ is better represented by viewing it as a linear code in $(\bbZ / 2 \bbZ)^k$.
Similarly, a subdirect product of $C_p^k$ leads to a linear code in $(\bbZ / p \bbZ)^k$. 
Thus, one has to relate the normaliser to objects studied in coding theory.

\subsection{Proposition}
Let $G \subset S_2^k (\subset S_{2k})$ be a subdirect product. 
Denote by $\phi \colon S_2^k \rightarrow (\bbZ / 2 \bbZ)^k$ the isomorphism that maps
the transposition $(2i-1,2i)$ to the $i$-th standard vector $e_i \in (\bbZ / 2 \bbZ)^k$. 
Let $A \subset S_k$ be the code automorphism group of $\phi(G)$.
Then the normaliser of $G$ in $S_{2k}$ is given by $S_2^k \rtimes A$. 

\proof This is part of~\cite[Thm.~5.1.15]{Ch}. \qed

\subsection{Remark}
More generally, for a prime $p > 2$ and a subdirect product $U \subset C_p^k$, 
the image of the normaliser of $U$ in $\Aut(C_p^n)$ is the automorphism group of the code $C$ 
corresponding to $U$ when we map $C_p^k$ to $(\bbZ / p \bbZ)^k$. Therefore, the normaliser
is given by $C_p^k \rtimes \Aut(C)$. More details are given in~\cite[Thm.~5.1.15]{Ch}.

\subsection{Remark}
One can work with subdirect products $G \subset (C_p \rtimes C_d)^k$ (for $d \mid (p-1)$) in a similar way.
First, one can use the image in $C_q^k$ for each prime $q \mid d$. Second, the $p$-Sylow subgroup of $G$ 
is characteristic. Therefore, the normaliser of $G$ is contained in the normaliser of the $p$-Sylow subgroup.
Thus, we get several normaliser overgroups and can form the intersection.

\subsection{Remark}
The above is not limited to groups having orbits of prime length. For example, denote by $K_4$ the Klein four-group 
(i.e., the regular representation of $(\bbZ / 2 \bbZ)^2$).  Now, let $U \subset K_4^n$ be a subdirect product.
Then one gets $S_4 \wr S_n$ as a first overgroup of the normaliser. Choosing $\pi_6 \colon S_4 \rightarrow S_6$
as the action on the cosets of a cyclic subgroup of order 4 results in degree 6 representations of $S_4$ and $K_4$. The
image of $K_4$ has 3 orbits of length 2. We denote by $\Pi_6 \colon S_4 \wr S_n \rightarrow S_6 \wr S_n$
the induced permutation representation. Then $\Pi_6(U)$ has only orbits of length 2 and therefore the
code approach can be used to compute the symmetric normaliser of $\Pi_6(U)$ and the normaliser of $U$ in $S_4 \wr S_n$ can 
easily be determined from this.
The remainder of this section turns this into a general algorithm.

\subsection{Algorithm\label{AffQuSearch}} (1-dimensional affine images) \\
{\bf Input:} Permutation groups $U \triangleleft N$ and a prime $p$. \\
{\bf Output:} A list of subgroups $H_i$ of $N$ such that the cost action of $U$ on $N/H_i$
results in a subdirect product of 1-dimensional affine groups over $\bbZ / p \bbZ$. \\
{\bf Steps:}
\begin{enumerate}
\item
Compute the index $p$ subgroups of $U$. 
Keep the subgroups $V \subset U$ such that the action of $U$ on $U / V$ 
is a 1-dimensional affine group.
\item \label{NClasses}
Partition the subgroups found into conjugacy classes with respect to conjugation by $N$.
Store a representative of each class.
\item \label{StepSearchH}
For each representative $R$ obtained in Step~\ref{NClasses}, find a subgroup $H \subset N$ of minimal index such that $H \cap U = R$.
Store the group $H$ found into a list.
\item
Return the list of subgroups $H$.
\end{enumerate}
 
\subsection{Remark}
As $U$ is normal in $G$, the group $R$ is normal in $H$. Therefore, it suffices in Step~\ref{StepSearchH} 
to search for $H$ in $N_N(R)$. If $V$ is not normal in $U$ then we have $H = N_G(R)$. Otherwise, $H$ can
easily be found by enumerating low index subgroups in $N_N(R) / R$. 

\subsection{Algorithm} (Overgroup refinement by code analysis) \label{GroupRefinedByCode} \\
{\bf Input:} A prime $p$, a permutation group $U \subset G_1^{e_1} \times \cdots \times G_m^{e_m}$,
 a normaliser overgroup $H$ contained in a group of the shape $W := N_1 \wr S_{e_1} \times \cdots \times N_m \wr S_{e_m}$
with the same orbits. Furthermore, the $N_i$ are assumed to be transitive and $G_i$ is normal in $N_i$.  \\
{\bf Output:} An overgroup of the symmetric normaliser of $U$. \\
{\bf Steps:}
\begin{enumerate}
\item
Call Algorithm~\ref{AffQuSearch} for the prime $p$ and each pair of 
groups $G_i \subset N_i$ in order to get a list of subgroups $R_{i,1},\ldots,R_{i,l_i}$, for each $i = 1,\ldots,m$.
\item
For each group $R_{i,j}$, use the coset action of $N_i$ on $N_i / R_{i,j}$ to induce a 
permutation representation of $\phi_{i,j}$ of $N_i \wr S_{e_i}$. 
Form the direct product $\Phi$ of the homomorphisms $\phi_{i,j} \circ \pi_i$.
Here, $\pi_i$ denotes the projection $W \rightarrow N_i \wr S_{e_i}$. 
\item
Compute the symmetric normalizer $S$ of the $p$-Sylow subgroup of $\Phi(U)$ using automorphisms
of linear codes.
\item
Return the intersection of $\Phi^{(-1)}(S \cap \Image(\Phi))$ with $H$  
as a refined overgroup of the normaliser.
\end{enumerate}

\subsection{Remark}
Homomorphisms can be used in even more ways. For example, if two homomorphisms $\phi \colon G \rightarrow Q$, 
$\psi \colon G \rightarrow R$, and a subgroup $U \subset G$ are given then we have
$$
N_G(U) \subset N_G(U \cap \ker(\psi) )
$$
and, more generally, 
$$
N_G(U) \subset \phi^{-1}(N_{Q}(\phi(U \cap \ker(\psi))))\, .
$$
In particular, when the number of factors in the map $\Phi$ in Algorithm~\ref{GroupRefinedByCode} is large,
the dimension of the code is large and the computation of the automorphism group may slow down. Then one can apply this 
concept and use some of the factors $\phi_{i,j}$ to determine the kernel and the remaining ones to map 
to a linear code of smaller dimension. In practice, the selection of the maps is done randomly and
the automorphism group is only computed for codes with at most $10^6$ words and for codes that
are images of a single factor $\phi_{i,j}$.

\section{Final algorithm and performance test}

\subsection{Introduction}
The author's implementation will be part of the next {\tt magma} version. It can be 
called via $\tt SymmetricNormalizer$ with the option \mbox{{\tt Al := "Chain"}}. 
The subalgorithms above  are called in it as described below.

\subsection{Algorithm\label{NorAlgTran}}(Normaliser for a transitive permutation group.)\\
{\bf Input:} A transitive permutation group $G \subset S_n$. \\
{\bf Output:} The normaliser of $G$ in $S_n$. \\
{\bf Steps:}
\begin{enumerate}
\item
If $G$ is primitive or $G$ is transitive of small degree (e.g. $n \leq 35$) or $\#G$ is small 
(e.g. $\leq 6 \cdot n$) then call the internal normaliser algorithm directly and return 
the result, as it is expected to be fast.
\item
Use  Algorithm~\ref{GraphGroup} to determine an overgroup $H$ of the normaliser.
\item
For each block system of $H$, denote by $\phi$ the action on the blocks and
replace $H$ by $\phi^{-1}(N_{\phi(H)}(\phi(G)))$.
\item
For each block system $(B_1,\ldots,B_k)$ of $H$, compute $S := \Stab_H(B_1)$.
If $\pi_{B_1}(G \cap S)$ is not normal in $\pi_{B_1}(S)$ then apply Proposition~\ref{F-Case} to refine $H$. 
\item
If the index of $G$ in $H$ is still large then use all the minimal  block systems to build
 block quotient maps $\phi$  (c.f.\ Remark~\ref{BlockQuotient}) and replace $H$ by 
$\phi^{-1}(N_{\phi(H)}(\phi(G)))$.
\item
Compute the normaliser of $G$ in $H$ with the internal algorithm and return it.
\end{enumerate}

\subsection{Algorithm\label{NorAlgIn}}(Normaliser for an intransitive permutation group.)\\
{\bf Input:} A intransitive permutation group $G \subset S_n$. \\
{\bf Output:} The normaliser of $G$ in $S_n$. \\
{\bf Steps:}
\begin{enumerate}
\item
If $G$ is of small degree (e.g. $n \leq 24$) then call the internal normaliser 
algorithm directly and return the result, as it is expected to be fast.
\item
Call Algorithm~\ref{DPWP} to determine an overgroup $H$ of the normaliser.
\item
Call Algorithm~\ref{GroupRefinedByCode} for each prime $p$ that divides the order of $G$
and replace $H$ by the overgroup obtained. 
\item
Replace $H$ by $\phi^{-1}(N_{\phi(H)}(\phi(G)))$, for all accessible, 
non-trivial homomorphisms $\phi$ with domain $H$.
The replacements of $H$ may result in new accessible, non-trivial homomorphisms. 
Run through the homomorphisms in the following order:
\begin{enumerate}
\item
Use the action on each orbit of the normaliser overgroup 
and apply the strategy described for transitive groups.
\item
Select a minimal system of blocks for each orbit with non-primitive action and use the
direct product of the actions.
\item
Map the normaliser overgroup to the action on its first $k$ orbits. Start with $k = 1$ and
increment $k$ until the normaliser of $G$ is obtained.
\end{enumerate}
\item
Return the computed normaliser of $G$. 
\end{enumerate}

\subsection{Tests using the database of transitive groups}
The database of transitive groups is available up to degree 48. 
In principle, one could use all the groups in the database for a test.
But, in degrees where there are more than 20000 groups, a random selection of groups was taken. 
The timings to compute the symmetric normalisers are shown in the table below. 
\medskip

\begin{tabular}{c|c|c|c|c|c}
Degree &  \#groups    & Avr.~time     & Max.~time     & Avr.~time      & Max.~time \\  
       &  tested      & internal      & internal      & Alg.~\ref{NorAlgTran} & Alg.~\ref{NorAlgTran} \\
\hline
 36 & 20000 & 0.021 sec. & 69.97 sec. & 0.011 sec. & 0.14 sec. \\ 
 40 & 20000 & 0.089 sec. & 85.85 sec. & 0.026 sec. & 0.15 sec. \\ 
 42 & 9491 &  0.047 sec. & 14.45 sec. & 0.017 sec. & 0.91 sec. \\ 
 44 & 2113 &  2.357 sec. & 1345 sec. & 0.029 sec. & 0.23 sec. \\ 
 45 & 10923 & 0.377 sec. & 596 sec. & 0.024 sec. & 1.06 sec. \\ 
 48 & 20000 & 0.520 sec. & 328.5 sec. & 0.047 sec. & 1.79 sec. \\ 

\end{tabular}

\subsection{Test with direct products}
Assume that an intransitive group is the direct product of its orbit images. Then 
Algorithm~\ref{DPWP} determines the normaliser exactly. 
For a test, the normaliser $N_{S_{60}}(G^5)$ was computed for all 301 transitive groups $G$ of degree 12. 
Here, the internal algorithm 
takes 399.91 seconds in total. Using Algorithm~\ref{NorAlgIn} took only 2.70 seconds.
Similarly, computing $N_{S_{60}}(G^4)$, for the 104 transitive groups $G$ of degree 15, 
takes in total  21.47 seconds  with  the internal algorithm,
whereas Algorithm~\ref{NorAlgIn} takes only 0.75 seconds.

\subsection{Test with subdirect products}
For each transitive group $G$ of degrees 4 to 12 in the database of transitive groups, we
form the direct product $P := G^n$.
By choosing random elements in $P$, we generate a subgroup $T$ and compute the normaliser.
Different numbers of generators were used to explore a variety of subgroups in $G^n$.
As the internal algorithm can be slow in these examples, it was only 
applied to small cases.
The timings observed are shown in the table below. 
\medskip

\begin{center}
\begin{tabular}{c|c|c|r|r|r|r|r}
$\deg(g)$ & $n$ & $\deg(T)$ &  \#groups    & Avr.~time     & Max.~time  & Avr.~time         & Max.~time \\  
          &     &           &  tested      & internal      & internal   & Alg.~\ref{NorAlgIn} & Alg.~\ref{NorAlgIn} \\
\hline  
     7 &      4 &     28 &     14 & 0.004 sec. & 0.010 sec.  & 0.004 sec. & 0.010 sec. \\ 
     8 &      4 &     32 &    204 & 0.014 sec. & 0.870 sec.  & 0.024 sec. & 0.110 sec. \\ 
     9 &      4 &     36 &     90 & 0.014 sec. & 0.280 sec.  & 0.016 sec. & 0.070 sec. \\ 
    10 &      4 &     40 &    128 & 0.035 sec. & 0.990 sec.  & 0.025 sec. & 0.200 sec. \\ 
    11 &      4 &     44 &     16 & 0.069 sec. & 0.710 sec.  & 0.004 sec. & 0.020 sec. \\ 
    12 &      4 &     48 &   1110 & 0.873  sec. & 503.5 sec.  & 0.035 sec. & 0.270 sec. \\[1mm] 
     4 &      6 &     24 &     22 & 0.010 sec. & 0.040 sec.   & 0.013 sec. & 0.020 sec. \\ 
     5 &      6 &     30 &     10 & 0.026 sec. & 0.140 sec.   & 0.010 sec. & 0.020 sec. \\ 
     6 &      6 &     36 &     56 & 0.049 sec. & 0.280 sec.   & 0.021 sec. & 0.050 sec. \\ 
     7 &      6 &     42 &     14 & 0.029 sec. & 0.130 sec.   & 0.021 sec. & 0.080 sec. \\ 
     8 &      6 &     48 &    346 & 1060 sec.  & 300950 sec.  & 0.055 sec. & 0.200 sec. \\ 
     9 &      6 &     54 &    134 & 27.07 sec. &  1912 sec.   & 0.046 sec. & 0.170 sec. \\ 
    10 &      6 &     60 &    198 & 15.91 sec. &  1843 sec.   & 0.075 sec. & 0.650 sec. \\ 
    11 &      6 &     66 &     16 & 1.509 sec. &  20.64 sec.  & 0.029 sec. & 0.180 sec. \\[1mm]
     4 &      8 &     32 &     36 & 2.336 sec. &  62.69 sec.  & 0.019 sec. & 0.04 sec. \\ 
     5 &      8 &     40 &     18 & 4.260 sec. &  54.49 sec.  & 0.022 sec. & 0.04 sec. \\ 
     6 &      8 &     48 &     92 & 2.031 sec. &  23.38 sec.  & 0.042 sec. & 0.19 sec. \\ 
     7 &      8 &     56 &     24 & 2.326 sec. &  21.74 sec.  & 0.031 sec. & 0.12 sec.  
\end{tabular}
\bigskip

\begin{tabular}{c|c|c|r|r|r}
$\deg(g)$ & $n$ & $\deg(T)$ &  \#groups     & Avr.~time         & Max.~time \\  
          &     &           &  tested       & Alg.~\ref{NorAlgIn} & Alg.~\ref{NorAlgIn} \\
\hline
     4 &     10 &     40 &     50       & 0.032 sec. & 0.080 sec. \\ 
     5 &     10 &     50 &     26       & 0.040 sec. & 0.090 sec. \\ 
     6 &     10 &     60 &    128       & 0.090 sec. & 0.490 sec. \\ 
     7 &     10 &     70 &     34       & 0.064 sec. & 0.260 sec. \\ 
     8 &     10 &     80 &    698       & 0.220 sec. & 1.39 sec. \\ 
     9 &     10 &     90 &    306       & 0.188 sec. & 1.25 sec. \\ 
    10 &     10 &    100 &    434       & 0.370 sec. & 4.86 sec. \\ 
    11 &     10 &    110 &     36       & 0.305 sec. & 1.73 sec. \\ 
    12 &     10 &    120 &   3948       & 0.879 sec. & 10.5 sec. \\[1mm]
     4 &     15 &     60 &     85       &  0.123 sec. & 0.340 sec. \\ 
     5 &     15 &     75 &     46       &  0.227 sec. & 0.800 sec. \\ 
     6 &     15 &     90 &    218       &  0.543 sec. & 3.190 sec. \\ 
     7 &     15 &    105 &     59       &  0.415 sec. & 2.410 sec. \\ 
     8 &     15 &    120 &   1138       &  2.041 sec. & 63.690 sec. 
\end{tabular}
\end{center}
\bigskip

As the groups were constructed randomly, repeated runs result in different timings.
The most surprising observation is that the average running time of the backtrack search 
is dominated by the slowest example, whereas the new approach results in more uniform timings.

\section{Conclusion}

The new normaliser algorithm performs well for a large variety of examples. 
Most notably, the deviation of the running time between similar examples is 
much smaller compared to the internal algorithm. Intransitive examples
show the most significant improvements.
This compares with the observations reported in~\cite{Ch}.
The use of graph automorphisms and the systematic usage of homomorphisms to groups, for which the code 
based approach applies, seem to be the most powerful new tools for the determination of
normalisers. The limit of accessible examples seems to be pushed at least to degree 100.

\subsection*{Further work}
A close inspection of the transitive groups, for which the new algorithm performs most slowly show the
following commonality. The graph approach of Algorithm~\ref{GraphGroup} results in a group with 
less block systems than the final normaliser. As a consequence, the overgroup is far too big and the 
refinement steps have too few block systems 
to use and, finally, the internal algorithm is called for a pair of groups of large index. 
Therefore, a search for more invariants that stop the normaliser from swapping block systems 
seems  promising. However, the examples found so far give no indication what better
invariants could be. Further, one could come up with more methods to generate accessible 
quotients in order to extend the chain of intermediate groups.
\medskip

Finally, it might be possible to use the approach for computing normalisers 
in other classes of groups, such as polycyclic groups or matrix groups.

\end{document}